\documentclass[12pt]{amsart}

\usepackage{amssymb}

\usepackage{enumerate}

\usepackage{graphicx}

\usepackage{multirow}

\makeatletter
\@namedef{subjclassname@2010}{%
  \textup{2010} Mathematics Subject Classification}
\makeatother

\newtheorem{theorem}{Theorem}[section]

\theoremstyle{definition}

\newtheorem{remark}[theorem]{Remark}

\numberwithin{equation}{section}

\DeclareMathOperator{\Sel}{Sel}
\DeclareMathOperator{\rank}{rank}
\DeclareMathOperator{\genus}{genus}

\DeclareMathOperator{\Jac}{Jac}

\newcommand{\Q}{{\mathbb{Q}}}
\newcommand{\Z}{{\mathbb{Z}}}
\newcommand{\bP}{{\mathbb{P}}}

\begin{document}


\baselineskip=17pt



\title{On arithmetic progressions on Edwards curves}
\author{Enrique Gonz\'alez--Jim\'enez}
\address{Universidad Aut{\'o}noma de Madrid, Departamento de Matem{\'a}ticas and Instituto de Ciencias Matem{\'a}ticas (ICMat), Madrid, Spain.}
\email{enrique.gonzalez.jimenez@uam.es}
\urladdr{http://www.uam.es/enrique.gonzalez.jimenez}

\maketitle

\begin{abstract}
 Assume $m\in\Z_{>0}$ and $a,q\in\Q$. Denote by $\mathcal{AP}_{m}(a,q)$ the set of rational numbers $d$ such that $a,a+q,\dots,a+(m-1)q$ form an arithmetic progression in the Edwards curve $E_d\,:\,x^2+y^2=1+d\,x^2 y^2$. In these conditions, we study the set $\mathcal{AP}_{m}(a,q)$ and we parametrize it by the rational points of an algebraic curve.
\end{abstract}

\section{Introduction}
Let $F(x,y)\in \Q[x,y]$ be a polynomial in two variables such that its locus defines a plane model of an elliptic curve $E$ over $\Q$. We say that a rational number $x$ belongs to $E(\Q)$ if $x$ is the $x$-coordinate of a point $P\in E(\Q)$. We also say that $P_1,\dots,P_n\in E(\Q)$ are in arithmetic progression if the corresponding $x$-coordinates $x_1,\dots,x_n$ form an arithmetic progression.  Several authors \cite{Bremner1999,Campbell2003,Garcia-Selfa-Tornero2005, Garcia-Selfa-Tornero2006, Alvarado2011,Schwartz-Solymosi-Zeeuw2011,Mohanty1975,Lee-Velez1992,Bremner-Silverman-Tzanakis2000,Ulas2012,Spearman2011,Ulas2005, MacLeod2006, Alvarado2010,Moody2011,Moody2011b} have studied this problem for different shapes of the polynomial $F(x,y)$, and some of them have worked with the $y$-coordinates instead of the $x$-coordinates. It is interesting to point that the shape of the polynomial $F(x,y)$ makes a big difference in this context. For example, if the polynomial $F(x,y)$ is symmetric in both variables then there is no difference between studying the points with respect to $x$-coordinates or to $y$-coordinates. This is the case of the so called Edwards curves, that is, when $F(x,y)=x^2+y^2-1-dx^2y^2$ for some $d\in\Q$, $d\neq 0,1$. In the sequel, such elliptic curves will be denoted by $E_d$. These curves have been deeply studied in cryptography and it has been found that the resulting addition formulas are very efficient, simple and symmetric (for instance, without distinction of addition and doubling).

In this paper we fix our attention on Edwards curves. The starting point of this work is Moody's paper \cite{Moody2011} where the case $0,\pm 1,\pm 2,\dots$ is studied, and in particular is proved that there are infinitely many choices of $d$ such that $0,\pm 1,\dots,\pm 4$ form an arithmetic progression in $E_d(\Q)$. At the end of his paper, Moody asked if this arithmetic progression could be longer and he tries with no success, by computer search, to find a rational $d$ fulfilling the extra requirement that $\pm 5$ belongs to the arithmetic progression too. We prove in this paper that a rational $d$ such that $0,\pm 1,\pm 2,\pm 3,\pm 4,\pm 5$ form an arithmetic progression in $E_d(\Q)$ does not exist. Moreover, Moody stated that it is an open problem to find an Edwards curve with an arithmetic progression of length $10$ or longer. Although we have found no answer to this question, we will try to convince the reader that the maximum possible length of an arithmetic progression in an Edwards curve is $9$.

Let $m\in\Z_{>0}$, $a, q\in\Q$ be such that $q> 0$, and denote
$$
\mathcal{AP}_{m}(a,q)=\{d\in\Q\,|\, \mbox{$a,a+q,a+2q,\dots,a+(m-1)q$ in $E_d(\Q)$}\}.
$$
Note that $\mathcal{AP}_{m}(a,q)=\mathcal{AP}_{m}(a+(m-1)q,-q)$, so we can assume without loss of generality that $q>0$.

Let us restrict for a moment to the case of symmetric progressions, i.e progressions such that if an element belongs to the sequence then its opposite does. There are two possibilities: else $a=0$ (central) or $a=\pm q/2$ (non-central). Note that if $0,q,\dots,m q$ belong to $E_d(\Q)$, then also $-q,\dots,-m q$ do. Therefore we denote
$$
\mathcal{S}_{c}{\mathcal{AP}}_{2m+1}(q)=\mathcal{AP}_{2m+1}(-mq,q).
$$
Similarly, if $q/2,3q/2,\dots, (2m-1)/2q$ belong to $E_d(\Q)$, then $-q/2,-3q/2,\dots$,\\ $- (2m-1)/2q$ do, and we denote by
$$
\mathcal{S}_{nc}{\mathcal{AP}}_{2m}(q)=\mathcal{AP}_{2m}(-(2m-1)q/2,q).
$$
Therefore if we denote by $\mathcal{SAP}_{m}$ the set of rationals $d$ such that a symmetric arithmetic progression of length $m$ belongs to $E_d(\Q)$, we have
$$
\mathcal{SAP}_{m}(q)=\left\{
\begin{array}{cl}
\mathcal{S}_{c}{\mathcal{AP}}_{m}(q) & \mbox{if $m$ is odd},\\
\mathcal{S}_{nc}{\mathcal{AP}}_{m}(q)& \mbox{if $m$ is even}.
\end{array}
\right.
$$
\begin{theorem}\emph{\bf (Non-Symmetric Case)}\label{teo_general}
Let $m\in\Z_{>0}$ and $a,q\in\Q$ be such that  $q>0$ and $(a,q)$ does not correspond to a symmetric arithmetic progression. Then
\begin{itemize}
\item[(i)] $\#\mathcal{AP}_{m}(a,q)=\infty$ if $m\le 3$, except for maybe a finite number of pairs $(a,q)$.
\item[(ii)] $\#\mathcal{AP}_{4}(a,q)=\infty$ if and only if $a+kq\in\{\pm 1\}$ for some $k\in\{0,1,2,3\}$.
\item[(iii)]  If $m\ge 5$, then $\#\mathcal{AP}_{m}(a,q)<\infty$ for any pair $(a,q)$.
\end{itemize}
\end{theorem}

\begin{theorem}\emph{\bf (Central Symmetric Case)} \label{teo_CS}
Let $m\in\Z_{>0}$ odd and $q\in\Q_{>0}$. Then:
\begin{itemize}
\item[(i)]  $\#\mathcal{S}_{c}{\mathcal{AP}}_{m}(q)=\infty$ if $m\le 7$.	
\item[(ii)]  $\#\mathcal{S}_{c}{\mathcal{AP}}_{9}(q)=\infty$ if and only if  $q\in\left\{1,\frac{1}{2},\frac{1}{3},\frac{1}{4}\right\}$.
\item[(iii)]  If $m\ge 11$ and $q\in\left\{1,\frac{1}{2},\frac{1}{3},\frac{1}{4}\right\}$, then $\#\mathcal{S}_{c}{\mathcal{AP}}_{m}(q)=0$.
\end{itemize}
\end{theorem}

\begin{theorem} \emph{\bf (Non-Central Symmetric Case)} \label{teo_NCS}
Let $m\in\Z_{>0}$ even and $q\in\Q_{>0}$. Then:
\begin{itemize}
\item[(i)]  $\#\mathcal{S}_{nc}{\mathcal{AP}}_{m}(q)=\infty$ if $m\le 6$.	
\item[(ii)]  $\#\mathcal{S}_{nc}{\mathcal{AP}}_{8}(q)=\infty$ if and only if $q\in\left\{2,\frac{2}{3},\frac{2}{5},\frac{2}{7}\right\}$.
\item[(iii)]  If $m\ge 10$ and $q\in\left\{2,\frac{2}{3},\frac{2}{5},\frac{2}{7}\right\}$, then $\#\mathcal{S}_{nc}{\mathcal{AP}}_{m}(q)=0$.
\end{itemize}
\end{theorem}

 A computer search was undertaken in Section \ref{sec_comp} to find a $q$ such that the set $\mathcal{SAP}_{m}(q)$ is non-empty for $m\ge 10$, but it was not succesful. So we left the following questions to the reader:\\

\noindent {\bf Question:\,}{\it Is $9$ the maximum length of an arithmetic progression on an Edwards curve, or in other words, is  $\#\mathcal{AP}_{m}(a,q)=0$ for any pair $a,q$ and $m\ge 10$?}

\section{Arithmetic-Algebraic-Geometric translation.}\label{translate}
Let $d\in\Q$ be such that $d\ne 0,1$. Then the Edwards curve is the elliptic curve defined by
$$
E_d:\,x^2+y^2=1+d\,x^2 y^2.
$$
We have that $(\pm 1,0),(0,\pm 1)\in E_d(\Q)$ (trivial points in the sequel). Moreover, since the model defined above is symmetric if  $(x,y)\in E_d(\Q)$, then $(\pm x, \pm y),(\pm y,\pm x)\in E_d(\Q)$.

If $(x,y)\in E_d(\Q)$ is a non-trivial point, then we can recover $d$ from $(x,y)$:
$$
d(x,y)=\frac{x^2+y^2-1}{x^2 y^2}.
$$
Assume that this point has the form $(x,y)=\left(a+nq,\frac{w}{z_n}\right)$, where $n\in\Z_{\ge 0}$, and $a,q\in \Q$ are such that $q\ne 0$. Then we define
$$
d_n:=d\left(a+nq,\frac{w}{z_n}\right)=\frac{w^2+z_n^2((a+nq)^2-1)}{(a+nq)^2w^2}.
$$
Notice that $a+nq\ne 0,\pm 1$ and $w\ne\pm z_n$ (resp. $n,q,w\ne 0$) since $d_n\ne 0,1$ (resp. the point is non-trivial).

Now, denote $\mathcal S=\{n_0,\dots,n_{m-1}\}\subset \Z_{\ge 0}$.  Then the finite set of equations
$$
\mathcal C_{\mathcal S}^{a,q}\, :\,\{d_i=d_j \,|\, i,j\in \mathcal S,\,i<j\}
$$
defines a curve in $\mathbb P^{m}$, where the points are $[w:z_0:\dots:z_{m-1}]$. Moreover, it is easy to check that a model for this curve may be obtained by fixing one element of $\mathcal S$, say $n_0$, and varying the rest of the elements of the set $\mathcal S$:
$$
\begin{array}{rl}
\mathcal C_{\mathcal S}^{a,q}\,:
\{(n_0 -n_j )q(2a+q(n_0+n_j)) w^2 + (a+n_jq)^2 (1 - (a+  n_0q)^2) z_{n_0}^2=&\\
 & \!\!\!\!\!\!\!\!\!\!\!\!\!\!\!\!\!\!\!\!\!\!\!\!\!\!\!\!\!\!\!\!\!\!\!\!\!\!\!\!\!\!\!\!\!\!\!\!\!\!\!\!\!\!\!\!\!\!\!\!\!\!\!\!\!\!\!\!\!\!\!\!\!\!\!\!\!\!\!\!\!\!\!\!\!\!\!\!\!\!\!\!\!\!\!\!\!\!\!\!\!(a+n_0q)^2 (1- (a+ n_jq)^2) z_{n_j}^2
\}_{j=1,\dots,m-1}.
\end{array}
$$
That is, $\mathcal C_{\mathcal S}^{a,q}$ is the intersection of $m-1$ quadric hypersurfaces in $\bP^{m}$ and therefore its genus is $(m-3)2^{m-2}+1$ (cf. \cite[Prop. 4]{Gonzalez-Jimenez-Xarles2013b} or \cite{Bombieri-Granville-Pintz1992}). Moreover, the points $[1:\pm1:\dots :\pm1]\in \mathcal C^{a,q}_{\mathcal S}$ correspond to the not-allowed case $d=1$. Therefore, if  $d\ne 0,1$, we obtain the following bijection:
$$
\mbox{\footnotesize
$
\left\{ \left(a+n_i q,\frac{w}{z_{n_i}}\right)\in E_d(\Q)\!\!\smallsetminus\!\!\{(\pm 1,0),(0,\pm 1)\}\,\Big|\,\, n_i\in \mathcal S \right\}\leftrightarrow  \mathcal C^{a,q}_{\mathcal S}(\Q)\smallsetminus \{[\pm1:\dots :\pm1]\}.$}
$$
\begin{remark}
Note that if $a$ and $q$ are not fixed, then $\mathcal C_{\mathcal S}^{a,q}$ is a variety of dimension $3$ in $\mathbb P^{m+2}$. Therefore, Edwards curves with $m$ points in arithmetic progressions are characterized by the rational points of a variety of dimension $3$. However, the computation of the whole set of rational points of a variety of dimension greater than one is still an intractable problem nowadays.
\end{remark}
We are going to rewrite the equations of $\mathcal C_{\mathcal S}^{a,q}$. For this purpose, and for any $i,j,k\in\Z_{> 0}$, we denote
$$
s_{i j} = \frac{q(i-j)(2a+(i+j)q)}{(a+iq)^2(1-(a+jq)^2)},\qquad r_{ij}=s_{ij}^{-1},\qquad t_{ijk}=\frac{s_{i k}}{s_{i j}}.
$$
Then
$$
\mathcal C_{\mathcal S}^{a,q}\,:\{ X_{j+1}^2= a_j X_0^2 + (1-a_j) X_1^2 \}_{j=1,\dots,m-1}
$$
where $a_j=s_{n_0n_j}$, $X_0=w$ and $X_{j+1}=z_{n_j}$ for any $n_j\in \mathcal S$.

Now we parametrize the first equation as
$$
[X_0:X_1:X_2]= [t^2-2t+a_1:-t^2+a_1:t^2-2a_1t+a_1].
$$
Using  this parametrization we substitute $X_0,X_1$ and $X_2$ in the rest of equations and we obtain a new system of equations of the curve, which depends on the parameter $t$:
$$
\mathcal C_{\mathcal S}^{a,q}\,:\{ X_{j+1}^2= t^4-4a_jt^3+2(-a_1+2a_j+2a_1a_j)t^2-4a_1a_jt+a_1^2 \}_{j=2,\dots,m-1}\,.
$$
Notice that each single equation defines an elliptic curve $\Q$-isomorphic to the elliptic curve with Weierstrass model
$$
\mathcal C^{a,q}_{\{n_0,n_1,n_j\}}\,:\,y^2=x(x+a_1-a_j)(x+a_j(a_1-1)).
$$
Here the isomorphism sends $[1:1:1:1]$ to $\mathcal O=[0:1:0]$, and if we denote by  $P_0=(0,0),P_1=(a_j-a_1,0),Q=(a_j,a_1a_j)$ then it sends the set $\{[\pm 1:\pm 1,\pm 1:\pm 1]\}$ to $\{ \mathcal O,P_1,P_2,P_1+P_2,Q,Q+P_1,Q+P_2,Q+P_1+P_2\}$.

Moreover, each pair of equations define a genus five curve $\mathcal C^{a,q}_{\{n_0,n_1,n_i,n_j\}}$ such that its jacobian $\Jac(\mathcal C^{a,q}_{\{n_0,n_1,n_i,n_j\}})$ splits completely over $\Q$ as the product of five elliptic curves. To check the previous assertion let us write $\mathcal C^{a,q}_{\{n_0,n_1,n_i,n_j\}}$ as (see \cite{Bremner1997}):
\begin{equation}\label{genus5}
\mathcal C^{a,q}_{\{n_0,n_1,n_i,n_j\}}\,:\,\left\{
\begin{array}{rcl}
X_2^2&\!\!\!\!=&\!\!\!\!b_2 X_0^2+(1-b_2)X_1^2,\\[1mm]
X_3^2&\!\!\!\!=&\!\!\!\!b_3X_0^2+(1-b_3)X_1^2,\\[1mm]
X_4^2&\!\!\!\!=&\!\!\!\!b_4 X_0^2+(1-b_4)X_1^2,
\end{array}
\right.
\end{equation}
where $X_3=X_i$, $X_4=X_j$ and $b_2=a_1$, $b_3=a_i$, $b_4=a_j$. Then we have five quotients of genus one such that each of them is the intersection of two quadric surfaces in $\bP^3$. Any elliptic curve $E_{(k)}$ has equations consisting of removing the variable $X_k$ from the previous system of equations. We display the Weierstrass model for those elliptic curves together with a (in general) non-torsion point on it:
$$
\begin{array}{ll}
E_{(4)}\,:\,y^2=x(x+b_2-b_3)(x+b_3(b_2-1)), & \!\!\!\!\!\!\!\!\!\!\!\!\!\!\!\!\!Q_4=(b_3,b_2b_3),\\
E_{(3)}\,:\,y^2=x(x+b_2-b_4)(x+b_4(b_2-1)), & \!\!\!\!\!\!\!\!\!\!\!\!\!\!\!\!\!Q_3=(b_4,b_2b_4),\\
E_{(2)}\,:\,y^2=x(x+b_3-b_4)(x+b_4(b_3-1)), & \!\!\!\!\!\!\!\!\!\!\!\!\!\!\!\!\!Q_2=(b_4,b_3b_4),\\
E_{(1)}\,:\,y^2=x\left(x+b_2(b_3-b_4)\right)\left(x+b_4(b_3-b_2)\right), & \!\! \!\!\!\!\!\!\!\!\!\!\!\!\!\!\!Q_1=(b_2 b_4,b_2 b_3b_4),\\
E_{(0)}\,:\,y^2=x\left(x+(b_2-1)(b_3-b_4)\right)\left(x+(b_4-1)(b_3-b_2)\right), & \\
&  \!\!\!\!\!\!\!\!\!\!\!\!\!\!\!\!\!\!\!\!\!\!\!\!\!\!\!\!\!\!\!\!\!\!\!\!\!\!\!\!\!\!\!\!\!\!\!\!\!\!\!\!\!\!\!\!\!\!\!\!\!\!\!\!\!\!\!\!\!\!\!\!\!\!\!\!\!\!\!\!\!\!\!\!\!\!\!\!\!\!Q_0=((b_4-1)(b_2-1),(b_2-1)(b_3-1)(b_4-1)).\\
\end{array}
$$
Therefore we have obtained $\Jac(\mathcal C^{a,q}_{\{n_0,n_1,n_i,n_j\}})\stackrel{\Q}{\sim} E_{(0)}\times E_{(1)}\dots \times E_{(4)}$ (cf. \cite{Bremner1997}). In general, $\rank_\Z\Jac(\mathcal C^{a,q}_{\{n_0,n_1,n_i,n_j\}})\ge 5=\genus(\mathcal C^{a,q}_{\{n_0,n_1,n_i,n_j\}})$, that is, the classical Chabauty's method \cite{Chabauty1941} does not work to obtain $C^{a,q}_{\{n_0,n_1,n_i,n_j\}}(\Q)$. The curve $\mathcal C^{a,q}_{\{n_0,n_1,n_i,n_j\}}$ has the same shape as the curve treated in \cite{Gonzalez-Jimenez-Xarles2013b} (with $m_0=b_2-1$, $m_1=-b_3$ and $m_2=-b_4$), where we developed a method based on covering collections and elliptic curve Chabauty techniques to obtain (under some hypotheses) the set of rational points of those curves (see \cite{Gonzalez-Jimenez2013} too). We intend now to apply this method to our curves, so let us write $\mathcal C^{a,q}_{\{n_0,n_1,n_i,n_j\}}$ in the following form:
$$
\mathcal C^{a,q}_{\{n_0,n_1,n_i,n_j\}}\,:\{ X_k^2= t^4-4b_kt^3+2(-b_2+2b_k+2b_2b_k)t^2-4b_2b_kt+b_2^2 \}_{k=3,4}.
$$
For $k\in\{3,4\}$ denote:
{\small
$$
\begin{array}{|c||c|c|l|}
\hline
l & \multicolumn{1}{|c|}{d_{k,l}}  &  \multicolumn{1}{|c|}{e_{k,l}} &  \multicolumn{1}{|c|}{p_{k,l,\pm}(t)} \\
\hline
1 & b_k(b_k-1) &  b_k(1-b_2) & t^2-2(b_k\pm\alpha_{k,1})t+b_2(-1+2(b_k\pm\alpha_{k,1}))\\
\hline
2 & (b_k-1)(b_k-b_2) & b_k-b_2 & t^2-2(b_k\pm\alpha_{k,2})t+b_2\\
\hline
3 & b_k(b_k-b_2)& 0 & t^2-2(b_k\pm\alpha_{k,3})t-b_2+2(b_k\pm\alpha_{k,3})\\
\hline
\end{array}
$$
}
where $\alpha_{k,l}=\sqrt{d_{k,l}}$. Next, choose $l_3,l_4\in \{1,2,3\}$ and for any $k\in\{3,4\}$ denote \\
$\bullet$ $\phi_k:E'_{(k)}\rightarrow E_{(k)}$ the $2$-isogeny
corresponding to the $2$-torsion point $(e_{k,l_k},0)\in E_{(k)}(\Q)$,\\
$\bullet$ $L=\Q(\alpha_{3,l_1},\alpha_{4,l_2})$,\\
$\bullet$ $\mathcal S_L(\phi_k)$ a set of representatives in $\Q$ of the image of the $\phi_k$-Selmer group $\Sel(\phi_k)$ in $L^*/(L^*)^2$ via the natural map,\\
$\bullet$  $\widetilde{\mathcal S_L}(\phi_3)$ a set of representatives of $\Sel(\phi_3)$ modulo the subgroup generated by the image of $[1:\pm 1:\pm 1:\pm 1:\pm 1]$ in this Selmer group,\\
$\bullet$  $\mathfrak{S} =\{\delta_3\delta_4 \ : \ \delta_3\in \widetilde{\mathcal S_L}(\phi_3), \delta_4\in \mathcal S_L(\phi_4)\}\subset \Q^*$,\\
$\bullet$  for any $\delta\in \mathfrak{S} $ and $s=(s_3,s_4)\in\{\pm\}\times \{\pm\}$ we define the genus one curve:
$$
H^\delta_s: \delta z^2=p_{3,l_3,s_3}(t)p_{4,l_4,s_4}(t).
$$
So we obtain the following:
$$
\begin{array}{ll}
\left\{t \in \bP^1(\Q)\, \Big|
\begin{array}{c}  \exists X_3,X_4\in \Q \mbox{ such that }\\
 (t,X_3,X_4)\in \mathcal C^{a,q}_{\{n_0,n_1,n_i,n_j\}}(\Q)
 \end{array}
 \right\} &
 \\[4mm]
 &\!\!\!\!\!\!\!\!\!\!\!\!\!\!\!\!\!\!\!\!\!\!\!\!\!\!\!\!\!\!\!\!\!\!\!\!\!\!\!\!\!\!\!\!\!\!\!\!\!\!\!\!\!\!\!\!\!\!\!\!\!\!\!\!\!\!\!\!\!\!\!\!\!\!\!
 \displaystyle \subseteq  \bigcup_{\delta \in \mathfrak{S}}
\left\{t \in \bP^1(\Q)\,\Big|
\begin{array}{l}
 \exists w \in L \mbox{ such that } (t,w) \in H_{s}^{\delta}(L)\\ \mbox{ for some $s\in\{\pm\}\times \{\pm\}$}
 \end{array}
 \right\}.
 \end{array}
$$
Note that in order to compute $\mathcal C^{a,q}_{\{n_0,n_1,n_i,n_j\}}(\Q)$ we must find a pair $l_3,l_4\in \{1,2,3\}$ such that for any $\delta \in \mathfrak{S}$ we can find $s\in \{\pm\}\times \{\pm\}$ where we can carry out all these computations to obtain the rational $t$-coordinates of $H_{s}^{\delta}(L)$. Before undertaking this task, however, we must face several problems, which in practice are solved by implementations in \verb|Magma| \cite{magma}: \\
$\bullet$  Is $H_s^{\delta}(L)$ empty? To answer this question we use the Bruin and Stoll's algorithm \cite{Bruin-Stoll2009}. If the answer is yes, we have finished with $\delta$ and go for another element of $\mathfrak S$. Otherwise, we must find (by brute force) a point on $H_s^{\delta}(L)$.\\
$\bullet$  Once we have found a point on $H_s^{\delta}(L)$, we use it to create an $L$-isomorphism with its Jacobian $\Jac(H_s^{\delta})$ and compute an upper bound for the rank $r$ of the Mordell-Weil group of the elliptic curve $\Jac(H_s^{\delta})(L)$. \\
$\bullet$  In the case where the rank $r< [L:\Q]$ we use the elliptic curve Chabauty algorithm (see \cite{Bruin2003}) to compute the $t$-coordinates of $H_{s}^{\delta}(L)$. For this purpose, we first must determine a system of generators of the Mordell-Weil group of $\Jac(H_s^{\delta})(L)$.

 If $\mathcal{S}=\{i,j,k,l\}$, then the curve $\mathcal C^{a,q}_{\mathcal{S}}$ has been defined by $\{d_i=d_j, d_i=d_k, d_i=d_l\}$. Note that we may describe this curve by $\{d_{n_1}=d_{n_2},d_{n_3}=d_{n_4}, d_{n_5}=d_{n_6}\}$ with $\{n_1,\dots,n_6\}=\{i,j,k,l\}$. If the order of the equations is assumed to be irrelevant, there are $16$ such descriptions; that is, we can consider $16$ models of $\mathcal C^{a,q}_{\mathcal{S}}$ of the form (\ref{genus5}). The possible values of $b_2,b_3,b_4$ (as a set) appear at table \ref{table_models16}. Then, we parametrized the first conic and make the appropriate substitutions on the other two conics. Therefore, if we take care now of the order of the equations, we have $48$ different models of $\mathcal C^{a,q}_{\mathcal{S}}$ of the form (\ref{genus5}). This is an important fact, since all the computations that we must carry out only may work (if they do) in a particular model. Notice that we only consider the case when $L$ is at most a quadratic field, as some of the computations are not well implemented for number fields of higher degree.

\begin{table}[h]
\begin{tabular}{|c|c||c|c||c|c||c|c|}
\hline
$\!N\!$ & $ \{b_2, b_3, b_4\} $ &$\!N \!$ & $ \{b_2, b_3, b_4\} $ &  $\!N\! $ & $ \{b_2, b_3, b_4\} $ & $N $ & $ \{b_2, b_3, b_4\} $  \\
\hline
\hline
$\!1\!$ &  $ s_{i j}  s_{i  k} ,   s_{i  l} $ &
$\!2\! $ & $ s_{j i} ,  s_{j  k} ,   s_{j l} $ &
$\!3\! $ & $ s_{k i} ,   s_{k j} ,   s_{k l} $ &
$\!4\!$ &  $ s_{l i} ,  s_{l j} ,   s_{l k}$ \\
\hline
$\!5\! $ &  $ r_{i  j} ,   t_{i j k} ,   t_{i j l} $ &
$\!6\! $ &  $ r_{i k}  ,   t_{i k j} ,     t_{i k l}$ &
$\!7 \!$ &  $ r_{i  l}  ,   t_{i l j} ,  t_{i l k} $ &
$\!8\! $ & $ r_{j  i}  ,     t_{j i k} ,     t_{j i l} $ \\
\hline
$\!9\! $ &  $ r_{j  k}  ,     t_{j k i}  ,     t_{j k l} $ &
$\!10\! $ & $ r_{j  l}  ,     t_{j l i}    ,     t_{j l k}$ &
$\!11\! $ & $ r_{k i}  ,     t_{k i j}    ,     t_{k i l}$ &
$\!12\! $ & $ r_{k j}  ,     t_{k j i}    ,     t_{k j l}$
\\ \hline
$\!13\! $ & $ r_{k l}  ,     t_{k l i}    ,     t_{k l j}$ &
$\!14\! $ & $ r_{l i}  ,     t_{l i j}    ,     t_{l i k}$ &
$\!15\! $ & $ r_{l j}  ,     t_{l j i}    ,     t_{l j k}$ &
$\!16\! $ & $ r_{l k}  ,     t_{l k  i}    ,     t_{l k  j}$
 \\
\hline
\end{tabular}
\caption{Models for $\mathcal{C}^{a,q}_{\{i,j,k,l\}}$}\label{table_models16}
\end{table}

We will use the following notation:
\begin{table}[h]
\begin{tabular}{|c|c|c|}
\hline
${\mathcal S}$  & $\mathcal C^{a,q}_{\mathcal S}$ &$ \genus(\mathcal C^{a,q}_{\mathcal S})$\\
\hline
\hline
 $\{i,j\}$   &      $\mathcal{C}_{ij}(a,q)$ & $0$\\
\hline
  $\{i,j,k\}$   &     $ \mathcal{E}_{ijk}(a,q)$ & $1$\\
\hline
              $\{i,j,k,l\}$   &      $\mathcal{D}_{ijkl}(a,q)$ & $5$\\
\hline
\hline
\end{tabular}
\caption{Notation for $\mathcal{C}^{a,q}_{\mathcal S}$}\label{tablanueva}
\end{table}

\section{Proof of Theorem \ref{teo_general}. Non-Symmetric case.}
We analyze in the sequel under which conditions a non-symmetric arithmetic progression $a,a+q,\dots, a+(m-1) q$ belongs to $E_d$. In particular $a\notin\{0,\pm q/2\}$. For this purpose, we are going to use the translation given in the previous section with $\mathcal S=\{0,1,\dots,m-1\}$. First notice that if $a+kq=0$ then this case corresponds to the central symmetric one. Now, if $a+kq\in\{\pm 1\}$ for some $k\in\mathcal S$, then we have $d_k=1$ and therefore we cannot use it and we should use instead the curve $\mathcal C^{a,q}_{\mathcal S^*}$; where ${\mathcal S}^*$ is the set $\mathcal S$ obtained by removing such values of $k$. Finally, if there exist $i,j\in\Z_{\ge 0}$ ,$i\ne j$, satisfying $a+iq=-1$ and $a+jq=1$, then the arithmetic progression must necessarily extend to a symmetric one. Therefore, we can assume that there is at most one value of $k\in\mathcal S$ satisfying $a+kq=1$ or $a+kq=-1$.

Let us prove the theorem depending on the set $\mathcal S^*$:

$\bullet$ $\#\mathcal S^*\le 1$: these cases are particularly simple. If $a=\pm 1$ then the set $\mathcal{AP}_m(a,q)$ is described by $d\ne 1$ when $m=1$ and by $d_1$ when $m=2$. Meanwhile, $d_0$ describes the case $m=1$ with $a\ne \pm 1$, and $m=2$ with $a+q=\pm 1$.

For the remaining cases (that is, when $\#\mathcal S^*> 1$), there is a bijection between the sets $\mathcal C^{a,q}_{\mathcal S^*}(\Q)$ and $\mathcal{AP}_m(a,q)$ for $m={\#\mathcal S}$.  Next table shows what is $\mathcal C^{a,q}_{\mathcal S^*}$ for each case (see Table \ref{tablanueva}).
\begin{table}[h]
\begin{tabular}{|c|c|c|c|c|c|}
\cline{2-6}
\multicolumn{1}{r}{} & \multicolumn{1}{|c|}{$m=1$} & \multicolumn{1}{|c|}{$m=2$} &\multicolumn{1}{|c|}{$m=3$} &\multicolumn{1}{|c|}{$m=4$} &\multicolumn{1}{|c|}{$m=5$}  \\
\hline
$a=\pm 1$ & $d\ne 1$ &  $d_1$ & $\mathcal{C}_{12}(a,q)$ & $\mathcal{E}_{123}(a,q)$ & $\mathcal{D}_{1234}(a,q)$\\
\hline
$a+q=\pm 1$ &  \multirow{4}{*}{$d_0$}  & $d_0$  & $\mathcal{C}_{02}(a,q)$ & $\mathcal{E}_{023}(a,q)$ & $\mathcal{D}_{0234}(a,q)$\\
\cline{1-1} \cline{3-6}
$a+2q=\pm 1$ &  &  \multirow{3}{*}{$\mathcal{C}_{01}(a,q)$} & $\mathcal{C}_{01}(a,q)$ & $\mathcal{E}_{013}(a,q)$ & $\mathcal{D}_{0134}(a,q)$\\
\cline{1-1} \cline{4-6}
$a+3q=\pm 1$ &  & &\multirow{2}{*}{$\mathcal{E}_{012}(a,q)$} & $\mathcal{E}_{012}(a,q)$ & $\mathcal{D}_{0124}(a,q)$\\
\cline{1-1} \cline{5-6}
 \multicolumn{1}{r|}{} & & & & $\mathcal{D}_{0123}(a,q)$  &\multicolumn{1}{r}{ } \\
\cline{2-5}\end{tabular}
\end{table}

We are going to split the proof depending on the cardinality of the set $\mathcal S^*$:\\

\noindent $\bullet$ $\mathcal S^*=\{i,j\}$: then the corresponding curve  is the conic $\mathcal{C}_{ij}(a,q)$ with equation
$$
\mathcal{C}_{ij}(a,q):z_j^2= s_{ij} w^2 + (1-s_{ij}) z_i^2.
$$
This conic has been parametrized on the previous section by
$$
[w:z_i:z_j]= [t^2-2t+s_{ij}:-t^2+s_{ij}:t^2-2s_{ij}t+s_{ij}],
$$
and therefore we have $\#\mathcal{AP}_m(a,q)=\infty$ when $m=\#\mathcal S$ and $\#\mathcal S^*=2$. These cases correspond to $\mathcal S=\{0,1\}$ and $a+kq\notin\{\pm 1\}$ for $k\in\{0,1\}$ or $\mathcal S=\{0,1,2\}$ and $a+kq\in\{\pm 1\}$ for $k\in\{0,1,2\}$.\\

\noindent $\bullet$ $\mathcal S^*=\{i,j,k\}$: we have proved on the previous section that the corresponding curve is an elliptic curve, i.e. it is $\Q$-isomorphic to the elliptic curve with Weierstrass model
$$
\mathcal{E}_{ijk}(a,q)\,:\,y^2=x(x+s_{ij}-s_{ik})(x+s_{ik}(s_{ij}-1)),
$$
and such that it has full $2$-torsion defined over $\Q$ and the extra rational point $Q=(s_{ik},s_{ij} s_{ik})$. Our first goal here is to prove that $Q$ is not a point of finite order for the cases
$$
(i,j,k,a)\in\{(1,2,3,\pm 1),(0,2,3,\pm 1-q),(0,1,3,\pm 1-2 q),(0,1,2,\pm 1-3q)\}.
$$
 By Mazur's theorem, $Q$ has infinite order if and only if $nQ$ is not a point of order $2$ for $n=1,2,3,4$, or equivalently, the $y$-coordinate $y_n$ of $nQ$ (that belongs to $\Q(q)$) is not $0$. We have factorized the numerator and denominator of $y_n$ for $n=1,2,3,4$ and obtained that the
factors of degree one correspond to symmetric arithmetic progressions. Therefore we have proved that $\mathcal{E}_{ijk}(a,q)$ has positive rank for any $(i,j,k,a)$ as above and any $q$ such that do not correspond to a symmetric arithmetic progression. Same arguments may be applied to the case $(i,j,k)=(0,1,2)$ and any $a,q$. In this case, $y_n\in \Q(a,q)$ and therefore the factors of its numerator and denominator define plane affine curves. All the corresponding genus zero curves come from the locus of the polynomials $a,q,a+q,2a+q,a+q\pm 1,a+2q\pm 1$. But these genus zero curves do not provide solutions since the possible rational points correspond to cases that have been excluded previously. The genus one curves define elliptic curves of rank zero and therefore only a finite number of points (in fact, the corresponding points are related to symmetric arithmetic progressions). The rest of the curves are of genus greater than one, and so they have only a finite number of rational points. In particular this concludes the proof of the first two statements of Theorem \ref{teo_general}.

To finish the proof of Theorem \ref{teo_general}, notice that if $m=\#\mathcal S\ge 5$ then $\#\mathcal S^*\ge 4$ and in this case the corresponding curve is of genus greater that one. Then, by Faltings' Theorem, this curve has a finite number of rational points. This proves that $\#\mathcal{AP}_m(a,q)<\infty$ when $m\ge 5$.

\section{Proof of Theorem \ref{teo_CS}. Central Symmetric case.}\label{sec_C}
Same arguments as above will be adapted to the central symmetric case. In this instance $a=0$, $\mathcal S=\{1,2,\dots,m\}$ and the condition $a+kq\in\{\pm 1\}$ becomes $kq=1$. Let ${\mathcal S}^*$ be the set $\mathcal S$ removing $k$.

If $\#\mathcal S^*\le 1$ the set  $\mathcal{S}_c\mathcal{AP}_{2s+1}(q)$ is described by the function $d_1$ when $s=1$ and $q\ne 1$; by $d\ne 1$ if $(s,q)=(1,1)$; by $d_2$ when $(s,q)=(2,1)$ and by $d_1$ when $(s,q)=(2,1/2)$.

If $\#\mathcal S^*\ge 2$, we use the bijection between $\mathcal C^{0,q}_{\mathcal S^*}(\Q)$ and $\mathcal{S}_c\mathcal{AP}_{2s+1}(q)$ for $s={\#\mathcal S}$.  Table \ref{table_moduli_central} shows what is $\mathcal C^{0,q}_{\mathcal S^*}$ for each case (see Table \ref{tablanueva}).\\
\begin{table}[h]
\begin{tabular}{|c|c|c|c|c|c|}
\hline
$q$ & \multicolumn{1}{|c|}{$m=3$} & \multicolumn{1}{|c|}{$m=5$} &\multicolumn{1}{|c|}{$m=7$} &\multicolumn{1}{|c|}{$m=9$} &\multicolumn{1}{|c|}{$m=11$}  \\
\hline
$1$ & $d\ne 1$ &  $d_2$ & $\mathcal{C}_{23}(0,1)$ & $\mathcal{E}_{234}(0,1)$ & $\mathcal{D}_{2345}(0,1)$\\
\hline
$1/2$ & \multirow{4}{*}{$d_1$}  & $d_1$ & $\mathcal{C}_{13}(0,1/2)$ & $\mathcal{E}_{134}(0,1/2)$ & $\mathcal{D}_{1345}(0,1/2)$\\
\cline{1-1} \cline{3-6}
$1/3$ &  & \multirow{3}{*}{$\mathcal{C}_{12}(0,q)$} & $\mathcal{C}_{12}(0,1/3)$ & $\mathcal{E}_{124}(0,1/3)$ & $\mathcal{D}_{1245}(0,1/3)$\\
\cline{1-1} \cline{4-6}
$1/4$ &  & & \multirow{2}{*}{$\mathcal{E}_{123}(0,q)$} & $\mathcal{E}_{123}(0,1/4)$ & $\mathcal{D}_{1235}(0,1/4)$\\
\cline{1-1} \cline{5-6}
 \multicolumn{1}{r|}{} & & & &  {$\mathcal{D}_{1234}(0,q)$} &\multicolumn{1}{r}{} \\
\cline{2-5}\end{tabular}
\caption{Moduli for $\mathcal{S}_c\mathcal{AP}_{m}(q)$}\label{table_moduli_central}
\end{table}

Now, if $\mathcal S^*=\{i,j\}$ the corresponding curve  is the conic $\mathcal{C}_{ij}(0,q)$ that has infinite number of points. Therefore we have $\#\mathcal{S}_c\mathcal{AP}_{2s+1}(q)=\infty$ when $s=\#\mathcal S$ and $\#\mathcal S^*=2$. These cases correspond to $\mathcal S=\{1,2,3\}$ and $q\in\{1,1/2,1/3\}$ or $\mathcal S=\{1,2\}$ and $q\notin\{1,1/2\}$.

The case $\mathcal S^*=\{i,j,k\}$ corresponds to the elliptic curve $\mathcal{E}_{ijk}(0,q)$ which has full $2$-torsion defined over $\Q$ and the extra rational point $Q=(s_{ik},s_{ij} s_{ik})$. Our objective is to prove that $Q$ is not a point of finite order for the cases $(i,j,k,q)\in\{(2,3,4,1),(1,3,4,1/2),(1,2,4,1/3)\}$ and $(i,j,k)=(1,2,3)$ for any $q\in\Q_{>0}$, $q\notin\{ 1,1/2,1/3\}$. The first attempt is to use the Nagell-Lutz theorem, so we compute an integral model of $\mathcal{E}_{ijk}(0,q)$ and we check if the coordinates of $Q'$ (the image  of the point $Q$ in this integral model) are not rational integers. The following table shows, for every case, an integral model and the $x$-coordinate of $nQ'$ for the first $n$ such that $nQ'$ has not integral coordinates.
$$
\begin{array}{|c|c|c|c|}
\hline
(i,j,k,q) &\mbox{integral model} & n & x(nQ')\\
 \hline
(2,3,4,1) & y^2 = x^3 - 25444800x - 35897472000 & 2 &185721/16 \\
 \hline
(1,3,4,1/2) & y^2 = x^3 - 11697075x + 15251172750 & 3 & 4532055/961 \\ \hline
(1,2,4,1/3) & y^2 = x^3 - 308700x - 55566000 & 1 & -4095/16 \\
\hline
\end{array}
$$
Therefore if $(i,j,k,q)\in\{(2,3,4,1),(1,3,4,1/2),(1,2,4,1/3)\}$ we infer that the point $Q$ is not of finite order.

Note that this procedure does not work for the case $(i,j,k)=(1,2,3)$ with $q\in\Q_{>0}$, $q\notin\{ 1,1/2,1/3\}$. By Mazur's theorem, $Q$ has infinite order if and only if $nQ$ is not a point of order $2$ for $n=1,2,3,4$ (that is, the $y$-coordinate $y_n$ of $nQ$, that belongs to $\Q(q)$, is not $0$. We factorized the numerator and denominator of $y_n$ for $n=1,2,3,4$ and we obtained that they have not any root $q\in\Q_{>0}$ with $q\notin\{ 1,1/2,1/3\}$. Thus, we have proved that $\mathcal{E}_{ijk}(0,q)$ has positive rank for any $(i,j,k,q)$ as above, and this proves that $\#\mathcal{S}_c\mathcal{AP}_{2s+1}(q)=\infty$ when $s=\#\mathcal S$ and $\#\mathcal S^*=3$. These cases correspond to $\mathcal S=\{1,2,3,4\}$ and $q\in\{1,1/2,1/3,1/4\}$ or $\mathcal S=\{1,2,3\}$ and $q\notin\{1,1/2,1/3,1/4\}$. In particular this concludes the proof of the statement: $\#\mathcal{S}_c\mathcal{AP}_{7}(q)=\infty$ for any $q\in\Q_{>0}$.  \\

The last case is $\mathcal S^*=\{i,j,k,l\}$, which corresponds to the genus five curve $\mathcal D_{ijkl}(0,q)$. By Falting's Theorem we have $\#\mathcal D_{ijkl}(0,q)(\Q)<\infty$. This proves that the set $\mathcal{S}_c\mathcal{AP}_{2s+1}(q)$ is finite when $s=\#\mathcal S$ and $\#\mathcal S^*=4$. These cases correspond to $\mathcal S=\{1,2,3,4,5\}$ and $q\in\{1,1/2,1/3,1/4\}$ or $\mathcal S=\{1,2,3,4\}$ and $q\notin\{1,1/2,1/3,1/4\}$. This concludes the proof of the fact that $\#\mathcal{S}_c\mathcal{AP}_{9}(q)=\infty$ if and only if  $q\in\{1,1/2,1/3,1/4\}$ and $\#\mathcal{S}_c\mathcal{AP}_{m}(q)<\infty$  for $m\ge 11$ and  any $q\in\Q_{>0}$.

\

In the remaining of this section, we will check that $\#\mathcal{S}_c\mathcal{AP}_{m}(q)=0$ if $q\in\{1,1/2,1/3,1/4\}$ for $m\ge 11$. Note that it is enough to prove it for $m=11$. That is, for $\mathcal S=\{1,2,3,4,5\}$ and $q\in\{1,1/2,1/3,1/4\}$ we have that the corresponding curve has genus five and therefore only a finite number of rational points. In fact, we are going to prove that $\mathcal C^{0,q}_{\mathcal{S}^*}(\Q)=\{[1:\pm1:\pm1:\pm1 :\pm1]\}$ for those values of $q$. For this purpose we apply the algorithm described on section \ref{translate}.

Let us start with the case $(a,q)=(0,1)$. Then the genus five curve is $\mathcal{D}_{2345}(0,1)$ and we choose the model $N=11$ from the table \ref{table_models16} with  $b_2=-4, b_3=7/32, b_4=-3/32$, and the pair $(l_3,l_4)=(1,2)$. In this case $L=\Q(\sqrt{-7})$, $\mathfrak{S}=\{\pm 1,\pm 10\}$ and the following polynomials:
$$
\begin{array}{l}
p_{3,1,+}(t)=t^2 + 1/16(-5\sqrt{-7} - 7)t + 1/4(-5\sqrt{-7} + 9),\\
p_{4,2,+}(t)= t^2 + 1/16(-25\sqrt{-7} + 3)t - 4.
\end{array}
$$
Now for any $\delta \in\mathfrak{S}$, we must compute all the points $(t,w) \in H_{\pm,\pm}^{\delta}(\Q(\sqrt{-7}))$ with $t\in \mathbb{P}^1(\Q)$ for some choice of the signs $s=(s_3,s_4)\in \{\pm\}\times\{\pm\}$ where
$$
H_s^{\delta}\,:\,\delta w^2=p_{3,1,s_3}(t)\,p_{4,2,s_4}(t).
$$
We have that $\rank_{\Z}H^{\pm 1}_{(+,+)}(\Q(\sqrt{-7}))=1$ therefore we can apply elliptic curve Chabauty to obtain the possible values of $t$. For $\delta=1$ (resp. $\delta=-1$) we obtain $t=\infty$ (resp. $t=-1$ ). For all of those values we obtain the trivial points $[1:\pm 1:\pm 1:\pm 1:\pm 1]\in \mathcal{D}_{2345}(0,1)(\Q)$. For $\delta\in\{\pm 10\}$, we obtain $H^{\delta}_{(+,+)}(\Q(\sqrt{-7}))=\varnothing$ using Bruin and Stoll's algorithm \cite{Bruin-Stoll2009}.

The following table shows all the previous data. Note that in the last column is stated if the corresponding points attached to $t$ in the curve are trivial or not:
$$
\begin{array}{c}
\begin{array}{|c|c|c|c|c|c|}
\hline
\delta & s & H^{\delta}_{s}(\Q(\sqrt{-7}))=\varnothing? & \rank_{\Z}H^{\delta}_{s}(\Q(\sqrt{-7})) & t & \mbox{trivial?}\\
\hline
1 & (+,+) & \mbox{no} & 1 & \infty & \mbox{yes}\\
\hline
-1 & (+,+) & \mbox{no} & 1 & 1 & \mbox{yes}\\
\hline
10 & (+,+) & \mbox{yes} & - & - & -\\
\hline
-10 & (+,+) & \mbox{yes} & - & - & -\\	
\hline
\end{array}\\[2mm]
\mbox{\footnotesize $(i,j,k,l)=(2,3,4,5), q=1, N=11, (b_2,b_3,b_4)=(-4,7/32,-3/32), (l_3,l_4)=(1,2)$}
\end{array}
$$
 From the previous table we obtain $\mathcal D_{2345}(0,1)(\Q)=\{[1:\pm 1:\pm 1:\pm 1:\pm 1]\}$; and therefore $\#\mathcal{S}_c\mathcal{AP}_{m}(1)=0$ for any $m\ge 11$. In particular, this answers one of Moody's questions\footnote{Moody asked \cite{Moody2011} if there exists $d\in\Q$, $d\ne 0,1$ such that  $0,\pm 1,\pm 2,\pm 3,\pm 4,\pm 5$ form an arithmetic progression in $E_d(\Q)$. Note, that after this paper was online (\texttt{http://arxiv.org/abs/1304.4361}), Bremner \cite{Bremner2013} has obtained a different proof of the non existence of such a $d$.}. 

The following three tables include the data related to the computation of all rational points of the curves
$\mathcal D_{ijkl}(0,q)(\Q)$ with $(i,j,k,l,q)\in\{(1,3,4,5,1/2)$,\\$(1,2,4,5,1/3),(1,2,3,5,1/4)\}$. In all these cases we have $\mathcal D_{ijkl}(0,q)=\{[1:\pm1:\pm1:\pm1 :\pm1]\}$.
$$
\begin{array}{c}
\begin{array}{c}
\begin{array}{|c|c|c|c|c|c|}
\hline
\delta & s & H^{\delta}_{s}(\Q(\sqrt{14}))=\varnothing? & \rank_{\Z}H^{\delta}_{s}(\Q(\sqrt{14})) & t & \mbox{trivial?}\\
\hline
1 & (+,-) & \mbox{no} & 1 & \infty& \mbox{yes}\\
\hline
2 & (+,+) & \mbox{yes} & - & - & -\\
\hline
-5 & (+,+) & \mbox{yes} & - & - & -\\
\hline
-10 & (+,+) & \mbox{yes} & - & - & -\\
\hline
\end{array}\\[2mm]
\mbox{\tiny $(i,j,k,l)=(1,3,4,5), q=1/2, N=3, (b_2,b_3,b_4)=(-3/25, 32/25, -64/125), (l_3,l_4)=(1,3)$}\\
\mbox{$\, $}
\\
\end{array}\\
\!\!\begin{array}{c}
\begin{array}{|c|c|c|c|c|c|}
\hline
\delta & s &\! H^{\delta}_{s}(\Q(\sqrt{21}))=\varnothing? \!&\! \rank_{\Z}H^{\delta}_{s}(\Q(\sqrt{21})) \!& \!t \!& \!\mbox{trivial?}\!\\
\hline
1 & (+,+) & \mbox{no} & 1 & \infty, \mbox{\small $27/25$} & \mbox{yes}\\
\hline
-1 & (+,+) & \mbox{yes} & - & - & -\\
\hline
6 & (+,+) & \mbox{yes} & - & - & -\\
\hline
-6 & (+,+) & \mbox{yes} & - & - & -\\
\hline
\end{array}\\[2mm]
\mbox{\tiny $(i,j,k,l)=(1,2,4,5), q=1/3, N=3, (b_2,b_3,b_4)=(27/25, 189/125, -81/175), (l_3,l_4)=(1,1)$}\\
\mbox{$\, $}
\\

\end{array}\\
\begin{array}{c}
\begin{array}{|c|c|c|c|c|c|}
\hline
\delta & s & H^{\delta}_{s}(\Q(\sqrt{105}))=\varnothing? & \rank_{\Z}H^{\delta}_{s}(\Q(\sqrt{105})) & t & \mbox{trivial?}\\
\hline
1 & (+,-) & \mbox{no} & 1 & \infty & \mbox{yes}\\
\hline
6 & (+,+) & \mbox{yes} & - & - & -\\
\hline
\end{array}\\[2mm]
\mbox{\tiny $(i,j,k,l)=(1,2,3,5), q=1/4, N=1, (b_2,b_3,b_4)=(128/3, -4, -128/7), (l_3,l_4)=(3,2)$}
\end{array}
\end{array}
$$
These computations conclude the proof of Theorem \ref{teo_CS}.

\section{Proof of Theorem \ref{teo_NCS}. Non-Central Symmetric case. }\label{sec_NC}

 We invoke the same arguments again. In this case we choose $a=-q/2$, $\mathcal S=\{1,2,\dots,m\}$ and the condition $a+kq\in\{\pm 1\}$ becomes $(2k-1)q=2$. Let be ${\mathcal S}^*$ the set $\mathcal S$ removing $k$.

If $\#\mathcal S^*\le 1$ the set  $\mathcal{S}_{nc}\mathcal{AP}_{m}(q)$ is described by the function $d_1$ when $m=2$ and $q\ne 2$; by $d\ne 1$ if $(m,q)=(2,2)$; by $d_2$ when $(m,q)=(4,2)$ and by $d_1$ when $(m,q)=(4,2/3)$.

Table \ref{table_moduli_non_central} identifies $\mathcal C^{-q/2,q}_{\mathcal S^*}$ when $\#\mathcal S^*>1$ (see Table \ref{tablanueva}).\\
{\small
\begin{table}[h]
\begin{tabular}{|c|c|c|c|c|c|}
\hline
$q$ & \multicolumn{1}{|c|}{\!$m=2$\!} & \multicolumn{1}{|c|}{\!$m=4$\!} &\multicolumn{1}{|c|}{$m=6$} &\multicolumn{1}{|c|}{$m=8$} &\multicolumn{1}{|c|}{$m=10$}  \\
\hline
$2$ & $d\ne 1$ &  $d_2$ & $\mathcal{C}_{23}(-1,2)$ & $\mathcal{E}_{234}(-1,2)$ & $\mathcal{D}_{2345}(-1,2)$\\
\hline
$2/3$ & \multirow{4}{*}{$d_1$}  & $d_1$ & $\mathcal{C}_{13}(-1/3,2/3)$ & $\mathcal{E}_{134}(-1/3,2/3)$ & $\mathcal{D}_{1345}(-1/3,2/3)$\\
\cline{1-1} \cline{3-6}
$2/5$ &  & \multirow{3}{*}{\!$\mathcal{C}_{12}(-q/2,q)$\!} & $\mathcal{C}_{12}(-1/5,2/5)$ & $\mathcal{E}_{124}(-1/5,2/5)$ & $\mathcal{D}_{1245}(-1/5,2/5)$\\
\cline{1-1} \cline{4-6}
$2/7$ &  & & \multirow{2}{*}{$\mathcal{E}_{123}(-q/2,q)$} & $\mathcal{E}_{123}(-1/7,2/7)$ & $\mathcal{D}_{1235}(-1/7,2/7)$\\
\cline{1-1} \cline{5-6}
 \multicolumn{1}{r|}{} & & & &  {$\mathcal{D}_{1234}(-q/2,q)$} &\multicolumn{1}{r}{} \\
\cline{2-5}\end{tabular}
\caption{Moduli for $\mathcal{S}_{nc}\mathcal{AP}_{m}(q)$}\label{table_moduli_non_central}
\end{table}
}

If $\#\mathcal S^*= 2$, then $\mathcal{S}_{nc}\mathcal{AP}_{2s}(q)$ is parametrized by a conic with infinitely many rational points. Therefore we have $\#\mathcal{S}_{nc}\mathcal{AP}_{2s}(q)=\infty$ when $s=\#\mathcal S$ and $\#\mathcal S^*=2$. These cases correspond to $\mathcal S=\{1,2,3\}$ and $q\in\{2,2/3,2/5\}$ or $\mathcal S=\{1,2\}$ and $q\notin\{2,2/3\}$.

Now, the elliptic curve $\mathcal{E}_{ijk}(-q/2,q)$ parametrizes the case when $\mathcal S^*=\{i,j,k\}$. This curve has all the $2$-torsion points defined over $\Q$ and the extra rational point $Q=(s_{ik},s_{ij} s_{ik})$. Using Nagell-Lutz we proved that $Q$ has infinite order for the cases
$$
(i,j,k,q)\in\{(2,3,4,2),(1,3,4,2/3),(1,2,4,2/5)\}.
$$
First we compute a suitable integral model of $\mathcal{E}_{ijk}(-q/2,q)$. Next table shows for every case the corresponding integral model and the $x$-coordinate of $nQ'$ for the first $n$ such that $nQ'$ has not integral coordinates (where $Q'$ is the image of $Q$ in this model):
$$
\footnotesize
\begin{array}{|c|c|c|c|}
\hline
\{i,j,k,q\} &\mbox{integral model} & n & x(nQ')\\
 \hline
\{2,3,4,2\} & y^2 = x^3 - 22427712x - 33269059584 & 3 & 2550847992/151321 \\
 \hline
\{1,3,4,2/3\} & y^2 = x^3 - 735300x + 242352000 & 2 & 18649/36 \\ \hline
\{1,2,4,2/5\} & y^2 = x^3 - 4615488 x - 3696371712 & 2 & 109761/25 \\
\hline
\end{array}
$$
The proof that $Q$ has infinite order in the case $(i,j,k,a,q)=(1,2,3,-q/2,q)$ with  $q\notin\{ 2,2/3,2/5\}$ is analogous to the case $(i,j,k,a,q)=(1,2,3,0,q)$ with  $q\notin\{ 1,1/2,1/3\}$ already discussed on the proof of Theorem \ref{teo_CS}. This proves that $\#\mathcal{S}_{nc}\mathcal{AP}_{2s}(q)=\infty$ when $s=\#\mathcal S$ and $\#\mathcal S^*=3$. These cases correspond to $\mathcal S=\{1,2,3,4\}$ and $q\in\{2,2/3,2/5,2/7\}$ or $\mathcal S=\{1,2,3\}$ and $q\notin\{2,2/3,2/5,2/7\}$. So we establish that $\#\mathcal{S}_{nc}\mathcal{AP}_{6}(q)=\infty$ for any $q\in\Q_{>0}$.  \\

Finally, the genus five curve $\mathcal D_{ijkl}(-q/2,q)$ corresponds to the case $\mathcal S^*=\{i,j,k,l\}$. Now, since $\#\mathcal D_{ijkl}(-q/2,q)(\Q)<\infty$ we obtain that the set $\mathcal{S}_{nc}\mathcal{AP}_{2s}(q)$ is finite when $s=\#\mathcal S$ and $\#\mathcal S^*=4$. These cases correspond to $\mathcal S=\{1,2,3,4,5\}$ and $q\in\{2,2/3,2/5,2/7\}$ or $\mathcal S=\{1,2,3,4\}$ and $q\notin\{2,2/3,2/5,2/7\}$. so we have proved that $\#\mathcal{S}_{nc}\mathcal{AP}_{8}(q)=\infty$ if and only if  $q\in\{2,2/3,2/5,2/7\}$ and $\#\mathcal{S}_{nc}\mathcal{AP}_{m}(q)<\infty$  for $m\ge 10$ and  any $q\in\Q_{>0}$.

The following four tables include the data related to the computation of all rational points of the curves $\mathcal D_{ijkl}(-q/2,q)$ with
$$
(i,j,k,l,q)\in\{(2,3,4,5,2),(1,3,4,5,2/3),(1,2,4,5,2/5),(1,2,3,5,2/7)\}.
$$
In all these cases we have $\mathcal D_{ijkl}(-q/2,q)(\Q)=\{[1:\pm1:\pm1:\pm1 :\pm1]\}$.
$$
\begin{array}{c}
\begin{array}{c}
\begin{array}{|c|c|c|c|c|c|}
\hline
\delta & s & H^{\delta}_{s}(\Q(\sqrt{15}))=\varnothing? & \rank_{\Z}H^{\delta}_{s}(\Q(\sqrt{15})) & t & \mbox{trivial?}\\
\hline
1 & (+,+) & \mbox{no} & 1 & \infty & \mbox{yes}\\
\hline
-1 & (+,+) & \mbox{yes} & - & - & -\\
\hline
6 & (+,+) & \mbox{no} & 1 & 1 & \mbox{yes}\\
\hline
-6 & (+,+) & \mbox{yes} & - & - & -\\	
\hline
\end{array}\\[2mm]
\mbox{\tiny $(i,j,k,l)=(2,3,4,5), a=-1,q=2, N=9, (b_2,b_3,b_4)=(7/5, 50, -4), (l_3,l_4)=(2,3)$}\\
\mbox{$\, $}
\\

\end{array}\\

\begin{array}{c}
\begin{array}{|c|c|c|c|c|c|}
\hline
\delta & s & H^{\delta}_{s}(\Q(\sqrt{10}))=\varnothing? & \rank_{\Z}H^{\delta}_{s}(\Q(\sqrt{10})) & t & \mbox{trivial?}\\
\hline
1 & (+,-) & \mbox{no} & 1 & \infty,0 & \mbox{yes}\\
\hline
-6 & (+,+) & \mbox{yes} & - & - & -\\	
\hline
\end{array}\\[2mm]
\mbox{\tiny $(i,j,k,l)=(1,3,4,5), a=-1/3,q=2/3, N=2, (b_2,b_3,b_4)=(7/25, 27/25, 27/125), (l_3,l_4)=(2,2)$}\\
\mbox{$\, $}
\\

\end{array}\\
\begin{array}{c}
\begin{array}{|c|c|c|c|c|c|}
\hline
\delta & s & H^{\delta}_{s}(\Q(\sqrt{21}))=\varnothing? & \rank_{\Z}H^{\delta}_{s}(\Q(\sqrt{21})) & t & \mbox{trivial?}\\
\hline
1 & (+,+) & \mbox{no} & 1 & \infty & \mbox{yes}\\
\hline
-1 & (+,+) & \mbox{yes} & - & - & -\\	
\hline
6 & (+,+) & \mbox{yes} & - & - & -\\	
\hline
-6 & (+,+) & \mbox{yes} & - & - & -\\	
\hline
\end{array}\\[2mm]
\mbox{\tiny $(i,j,k,l)=(1,2,4,5), a=-1/5,q=2/5, N=6, (b_2,b_3,b_4)=(5/7, 1/50, -1/4), (l_3,l_4)=(2,3)$}\\
\mbox{$\, $}
\\

\end{array}\\
\begin{array}{c}
\begin{array}{|c|c|c|c|c|c|}
\hline
\delta & s & H^{\delta}_{s}(\Q(\sqrt{7}))=\varnothing? & \rank_{\Z}H^{\delta}_{s}(\Q(\sqrt{7})) & t & \mbox{trivial?}\\
\hline
1 & (-,+) & \mbox{no} & 1 & \infty,0 & \mbox{yes}\\
\hline
2 & (+,+) & \mbox{yes} & - & - & -\\	
\hline
5 & (+,+) & \mbox{yes} & - & - & -\\	
\hline
10 & (+,+) & \mbox{yes} & - & - & -\\	
\hline
\end{array}\\[2mm]
\mbox{\tiny $(i,j,k,l)=(1,2,3,5), a=-1/7,q=2/7, N=1, (b_2,b_3,b_4)=(245/2, -49/5, -49), (l_3,l_4)=(2,2)$}
\end{array}
\end{array}
$$

This concludes the proof of theorem \ref{teo_NCS}.

\section{Some computations}\label{sec_comp}

We would like to find an arithmetic progression on an Edwards curve as large as possible. As in the symmetric ones fewer restrictions appear, we have undertaken a computer search on \texttt{Magma} to find a non-trivial rational point $P$ of height $H(P)\le 10^6$ in the curve $\mathcal D_{1234}(0,q)$ or in the curve $\mathcal D_{1234}(-q/2,q)$ for positive rationals $q$ of height $H(q)\le 100$ and $q\notin\{1,1/2,1/3,1/4\}$ or $q\notin\{2,2/3,2/5,2/7\}$ respectively. There are $6087$ such $q$'s. We used the following models for $\mathcal D_{1234}(0,q)$ and $\mathcal D_{1234}(-q/2,q)$:
$$
\mathcal D_{1234}(0,q)\,:\,\left\{
\begin{array}{lcl}
3X_0^2+4(q^2 - 1)X_1^2+(1-4q^2)X_2^2&=&0,\\[1mm]
8 X_0^2+9(q^2 - 1)X_1^2+(1-9q^2)X_3^2&=&0,\\[1mm]
15 X_0^2+16(q^2 - 1)X_1^2+(1-16q^2 )X_4^2&=&0,\\[1mm]
\end{array}
\right.
$$

$$
\mathcal D_{1234}(-q/2,q)\,:\,\left\{
\begin{array}{lcl}
32X_0^2+9(q^2 - 4)X_1^2+(4-9q^2)X_2^2&=&0,\\[1mm]
96 X_0^2+25(q^2 - 4)X_1^2+(4-25q^2)X_3^2&=&0,\\[1mm]
192 X_0^2+49(q^2 - 4)X_1^2+(4-49q^2 )X_4^2&=&0.\\[1mm]
\end{array}
\right.
$$
 We havenot found such a rational point. On the other hand, using the techinques of the proof of the last item of Theorems \ref{teo_CS} and \ref{teo_NCS}, we are able to prove that $\#\mathcal D_{1234}(0,q)=16$ for
$$q\in\left\{
\begin{array}{c}
19/11,11/13,49/46,13/3,3/2,3/7,2,11/43,1/11,\\
7/11,1/8,1/7,1/6,8/17,1/5,11/38,5/17,2/3,11/37,\\
7/13,59/61,29/53,3/4,11/19,3/8,37/95,11/28
\end{array}
\right\},
$$
and $\#\mathcal D_{1234}(-q/2,q)=16$ for
$$q\in\left\{\\
\begin{array}{c}
2/9,22/13,14,22/7,14/11,2/35,6/7,22/25,34/19,\\
2/17,2/15,22/73,62/33,2/13,38/35,10/7,34/49,22/31,\\
26/21,10/23,34/77,14/19,26/11,38/77,22/43,6/11
\end{array}\right\},
$$
Then for the corresponding list we have proved that
$$
\#\mathcal{S}_c\mathcal{AP}_{9}(q)=0\quad\mbox{and}\quad\#\mathcal{S}_{nc}\mathcal{AP}_{8}(q)=0$$
respectively.


\subsection*{Acknowledgements}
We would like to thank to Nils Bruin, Luis Dieulefait and Xevi Guitart for some interesting and useful discussions, and Jos\'e M. Tornero, who read the earlier versions of this paper carefully.

\

This research was partly supported by the grant MTM2012--35849.

\end{document}